\documentclass[11pt]{article}
\usepackage{tikz}
\usepackage[all]{xy}
\usepackage{amsfonts}
\usepackage{amsfonts}
\usepackage{amsfonts}
\usepackage{mathrsfs}
\usepackage{bm}
\usepackage{cite}
\usepackage[colorlinks,linkcolor=blue,citecolor=blue]{hyperref}
\usepackage{amssymb,amsmath} % font used for R in Real numbers
\usepackage{makecell}
\usepackage{tikz}
 \allowdisplaybreaks
\usetikzlibrary{arrows,shapes,chains}
 \allowdisplaybreaks

\catcode`!=11
\let\!int\int \def\int{\displaystyle\!int}
\let\!lim\lim \def\lim{\displaystyle\!lim}
\let\!sum\sum \def\sum{\displaystyle\!sum}
\let\!sup\sup \def\sup{\displaystyle\!sup}
\let\!inf\inf \def\inf{\displaystyle\!inf}
\let\!cap\cap \def\cap{\displaystyle\!cap}
\let\!max\max \def\max{\displaystyle\!max}
\let\!min\min \def\min{\displaystyle\!min}
\let\!frac\frac \def\frac{\displaystyle\!frac}
\catcode`!=12

\let\oldsection\section
\renewcommand\section{\setcounter{equation}{0}\oldsection}

\allowdisplaybreaks
\def\pf{\noindent \it{Proof.}\rm\quad}

\def\z{\zeta}

\newtheorem{thm}{Theorem}[section]
\newtheorem{lem}[thm]{Lemma}
\newtheorem{cor}[thm]{Corollary}

\newtheorem{pro}[thm]{Proposition}
\begin{document}
%%%%%%%%%%%%%%%%%%%% title %%%%%%%%%%%%%%%%%%%%%%%%%%%%%%%%%%%%%%%%%%%%%%%%
\title {\bf Weighted sum formulas of multiple $t$-values with even arguments}
\author{
{ Zhonghua Li$^{a,}$\thanks{Email:  zhonghua\_li@tongji.edu.cn}\quad Ce Xu$^{b,c,}$\thanks{Email: 19020170155420@stu.xmu.edu.cn; 9ma18001g@math.kyushu-u.ac.jp}}\\[1mm]
\small a. School of Mathematical Sciences, Tongji University\\
 \small Shanghai 200092, P.R. China\\
\small b. Multiple Zeta Research Center, Kyushu University \\
\small  Motooka, Nishi-ku, Fukuoka 819-0389, Japan\\
\small c. School of Mathematical Sciences, Xiamen University\\
\small Xiamen
361005, P.R. China}

\date{}
\maketitle \noindent{\bf Abstract} In this paper, we study the weighted sums of multiple $t$-values and of multiple $t$-star values at even arguments. Some general weighted sum formulas are given, where the weight coefficients are given by (symmetric) polynomials of the arguments.
\\[2mm]


\section{Introduction}

We begin with some basic notation. A finite sequence ${\bf k} = (k_1,\ldots, k_n)$ of positive integers is called an index. We
put
\[|{\bf k}|:=k_1+\cdots+k_n,\quad d({\bf k}):=n,\]
and call them the weight and the depth of ${\bf k}$, respectively. If $k_1>1$, ${\bf k}$ is called admissible. Let $I(k,n)$ be the set of all indices of weight $k$ and depth $n$.

For an admissible index ${\bf k}=(k_1,\ldots,k_n)$, the multiple zeta value and the multiple zeta-star value are defined by
\begin{align*}
\zeta({\bf k})\equiv\zeta(k_1,k_2,\ldots,k_n):=\sum\limits_{m_1>m_2>\cdots>m_n>0 } \frac{1}{m_1^{k_1}m_2^{k_2}\cdots m_n^{k_n}}
\end{align*}
and
\begin{align*}
\zeta^\star({\bf k})\equiv\zeta^\star(k_1,k_2,\ldots,k_n):=\sum\limits_{m_1\geqslant m_2\geqslant \cdots\geqslant m_n\geqslant 1 } \frac{1}{m_1^{k_1}m_2^{k_2}\cdots m_n^{k_n}},
\end{align*}
respectively. The systematic study of multiple zeta values began in the early 1990s with the works of Hoffman \cite{H1992} and Zagier \cite{DZ1994}. After that it has been attracted a lot of research on them in the last three decades (see, for example, the book of Zhao \cite{Z2016}).

Let $m,k,n$  be positive integers with $k\geqslant n$, and let $f(x_1,\ldots,x_n)\in \mathbb{Q}[x_1,\ldots,x_n]$ be a symmetric polynomial. Set
\[E_f(2m,k,n):=\sum\limits_{(k_1,k_2,\ldots,k_n)\in I(k,n)}f(k_1,k_2,\ldots,k_n)\zeta(2mk_1,2mk_2,\ldots,2mk_n),\]
which is a weighted sum of multiple zeta values with even arguments of weight $2mk$ and depth $n$. If $f(x_1,\ldots,x_n)=1$, we denote $E_f(2m,k,n)$ by $E(2m,k,n)$. The evaluations of these weighted sums $E_f(2m,k,n)$ have been attracted the attention of many researches. In \cite{GKZ2006}, Gangl, Kaneko and Zagier proved that $E(2,k,2)=\frac{3}{4} \z(2k)$. Later, Nakamura gave a different proof of this result in \cite{N2009}. Shen and Cai \cite{SC2012} studied the sums $E(2,k,3)$ and $E(2,k,4)$, and evaluated them in terms of $\z(2k)$ and $\z(2)\z(2k-2)$.  Using different methods, Hoffman \cite{H2017} and Gen${\check{\rm c}}$ev \cite{G2016} gave the explicit formula of $E(2,k,n)$. Furthermore,  Gen${\check{\rm c}}$ev \cite{G2016} proposed a conjecture on the weighted sum $E(4,k,n)$. More general sums $E(2m,k,n)$ with $m\geqslant 2$ have been considered by Komori, Matsumoto and Tsumura \cite{KMT2014}. And the explicit evaluation formula of $E(2m,k,n)$ was obtained just recently in \cite{ELO2017,ELO2018} and \cite{LQ2016}. Later in \cite{GL2015}, Guo, Lei and Zhao considered the weighted sums $E_f(2,k,2)$ and $E_f(2,k,3)$, and found that they can be evaluated by zeta values at even arguments. Moreover, they conjectured that
\begin{align*}
E_f(2,k,n)=\sum\limits_{l=0}^T c_{f,l}(k)\z(2l)\z(2k-2l),
\end{align*}
where $T={\rm max}\{[r+n-2]/2,[(n-1)/2]\}$, $c_{f,l}(x)\in \mathbb{Q}[x]$ depends only on $l$ and $f$, and with $\deg c_{f,l}(x)\leqslant \deg_{x_1}f(x_1,\ldots,x_n)$. Here $r=\deg f(x_1,\ldots,x_n)$ and for a real number $\alpha$, we denote by $[\alpha]$ the greatest integer that not exceeding $\alpha$. Recently, this conjecture was proved by the first author and Qin in \cite{LQ2019} with restriction $\deg c_{f,l}(x)\leqslant r+n-2l-1$. A similar weighted sum formula of the multiple zeta-star values with even arguments was simultaneously obtained by the first author and Qin in \cite{LQ2019}.

In a recent paper \cite{H2016}, Hoffman introduced and studied an odd variant of multiple zeta values,  which is defined for an admissible index ${\bf k}=(k_1,k_2,\ldots,k_n)$ as
\begin{align*}
t({\bf k})\equiv t(k_1,k_2,\ldots,k_n):=\sum\limits_{m_1>m_2>\cdots>m_n>0\atop m_i:\text{odd}} \frac{1}{m_1^{k_1}m_2^{k_2}\cdots m_n^{k_n}},
\end{align*}
and is called a multiple $t$-zeta value. Similarly, one can define a multiple $t$-star value by
\begin{align*}
t^\star({\bf k})\equiv t^\star(k_1,k_2,\ldots,k_n):=\sum\limits_{m_1\geqslant m_2\geqslant\cdots\geqslant m_n\geqslant 1\atop m_i:\text{odd}} \frac{1}{m_1^{k_1}m_2^{k_2}\cdots m_n^{k_n}}.
\end{align*}
Then similar as multiple zeta values, for any positive integers $m,k,n$ with $k\geqslant n$ and any symmetric
polynomial $f(x_1,\ldots,x_n)\in\mathbb{Q}[x_1,\ldots,x_n]$, we define the weighted sums of multiple $t$-values and of multiple $t$-star values by
\[T_f(2m,k,n):=\sum\limits_{(k_1,k_2,\ldots,k_n)\in I(k,n)}f(k_1,k_2,\ldots,k_n)t(2mk_1,2mk_2,\ldots,2mk_n)\]
and
\[T_f^\star(2m,k,n):=\sum\limits_{(k_1,k_2,\ldots,k_n)\in I(k,n)}f(k_1,k_2,\ldots,k_n)t^\star(2mk_1,2mk_2,\ldots,2mk_n).\]
If $f(x_1,\ldots,x_n)=1$, we set
\[T(2m,k,n)=T_f(2m,k,n),\qquad T^{\star}(2m,k,n)=T^\star_f (2m,k,n).\]
There are some work on the evaluations of the sums $T(2m,k,n)$. For example, using similar but more complicated ideas from \cite{SC2012} Shen and Cai gave a few sum formulas of $T(2,k,n)$ for $n\leqslant 5$ in \cite{SC2011}.
In \cite{Z2015}, Zhao gave two explicit formulas of $T(2,k,n)$. Furthermore, Shen and Jia \cite{SJ2017} gave some explicit evaluation formulas of $T(2m,k,n)$. We remark that as for multiple zeta values, one can obtain the evaluation formulas of $T(2m,k,n)$ and $T^{\star}(2m,k,n)$ algebraically. In fact, using \cite[Theorem 2.3]{H2016} and \cite[Proposition 3.26]{LQ2016}, one can express $T(2m,k,n)$ and $T^{\star}(2m,k,n)$ in terms of $t(2m,\ldots,2m)$ and $t^{\star}(2m,\ldots,2m)$. Then one gets evaluation formulas of $T(2m,k,n)$ and $T^{\star}(2m,k,n)$ from that of multiple $t$ and $t$-star values with all arguments equal to the same even number $2m$.

In this paper, using a similar method as in \cite{LQ2019}, we study the weighted sums $T_f(2,k,n)$ and $T_f^\star(2,k,n)$.  Our main result is the following theorem.

\begin{thm}\label{thm2}
Let $n,k$ be positive integers with $k\geqslant n$. Let $f(x_1,\ldots,x_n)\in \mathbb{Q}[x_1,\ldots,x_n]$ be a symmetric polynomial of degree $r$. Then we have
\begin{align}\label{1.2}
T_f(2,k,n)=\sum\limits_{l=0}^{{\rm min}\{T,k\}} c_{f,l}(k)\z(2l)t(2k-2l)
\end{align}
and
\begin{align}\label{1.3}
T_f^\star(2,k,n)=\sum\limits_{l=0}^{{\rm min}\{T,k\}} c_{f,l}^\star(k)\z(2l)t(2k-2l),
\end{align}
where $T={\rm max}\{[(r+n-2)/2],[(n-1)/2]\}$, $c_{f,l}(x),c_{f,l}^\star(x)\in \mathbb{Q}[x]$ depend only on $l$ and $f$, and with $\deg c_{f,l}(x), \deg c_{f,l}^\star(x)\leqslant r+n-2l-1$.
\end{thm}

To prove Theorem \ref{thm2}, we use the symmetric sum formulas of multiple $t$-values and of multiple $t$-star values \cite[Theorems 2.5 and 2.8]{H2016}. Then we find it is sufficient to study the weighted sums of the products of $t$-values at even integers and prove the following theorem.

\begin{thm}\label{thm1}
Let $n,k$ be positive integers with $k\geqslant n$. Let $f(x_1,\ldots,x_n)\in\mathbb{Q}[x_1,\ldots,x_n]$ be a polynomial of degree $r$. Then we have
\begin{align}\label{1.1}
\sum\limits_{k_1+\cdots+k_n=k\atop k_j\geqslant 1}f(k_1,\ldots,k_n)t(2k_1)\cdots t(2k_n)=\sum\limits_{l=0}^{\min\{T,k\}}e_{f,l}(k)\zeta(2l)t(2k-2l),
\end{align}
where $T=\max\{[(r+n-2)/2],[(n-1)/2]\}$, $e_{f,l}(x)\in\mathbb{Q}[x]$ depends only on $l$ and $f$, and with $\deg e_{f,l}(x)\leqslant r+n-2l-1$.
\end{thm}

Note that the polynomial $f(x_1,\ldots,x_n)$ in Theorem \ref{thm1} is not necessarily symmetric.

To prove Theorem \ref{thm1}, we use Euler's formula \cite[Eq. (1.4)]{H2016}, which expresses $t(2k)$ by the Bernoulli numbers. Then it is enough to treat the weighted sums of products of the Bernoulli numbers. And we do this by using the generating function of the Bernoulli numbers. We give the proofs of Theorem \ref{thm1} and Theorem \ref{thm2} in Section \ref{Sec:Proof}.

Although the weighted sum formulas \eqref{1.2}-\eqref{1.1} are not so concrete, one can get the explicit formulas for given positive integers $n,k$ and a given polynomial $f$ according to the procedure of our proof. We list some weighted sum formulas as examples in Appendix \ref{Sec:Example}.

\section{Proofs}\label{Sec:Proof}

\subsection{Preliminary knowledge}

We begin with the definition of the Bernoulli numbers. The generating function of the Bernoulli numbers $\{B_i\}$ is
$$\sum\limits_{i=0}^{\infty} \frac{B_i}{i!}x^i=\frac{x}{e^x-1}.$$
Let $\beta_i:=(2^i-1)B_i$ and
\begin{align*}
&F(x):=\frac{x}{2}-\frac{x}{e^x+1}.
\end{align*}
Then since $B_0=1$, $B_1=-\frac{1}{2}$ and $B_i=0$ for odd $i\geqslant 3$, we find that
\begin{align*}
&F(x)=\sum\limits_{i=1}^\infty\frac{\beta_{2i}}{(2i)!}x^{2i}=\sum\limits_{i=0}^\infty\frac{\beta_{2i}}{(2i)!}x^{2i}.
\end{align*}
Hence, $F(x)$ is an even function.

Let $D:=x\frac{d}{dx}$ and $H(x):=\frac{x}{e^x+1}$. Then we have
 $$D H(x)=(1-x)H(x)+H(x)^2.$$
Therefore, one can get the following theorem without difficulty.

\begin{thm}\label{thm2.1}
For any nonnegative integer $m$,
\begin{align}
D^mF(x)=\sum\limits_{i=0}^{m+1}F_{mi}(x)H(x)^i.
\label{Eq:Diff-f}
\end{align}
Here $F_{mi}(x)$ are polynomials determined by $F_{00}(x)=\frac{x}{2}$, $F_{01}(x)=-1$ and the recurrence relations
\begin{align}
\begin{cases}
F_{m0}(x)=xF_{m-1,0}'(x) & \text{for\;} m\geqslant 1,\\
F_{m,m+1}(x)=mF_{m-1,m}(x) & \text{for\;} m\geqslant 1,\\
F_{mi}(x)=xF'_{m-1,i}(x)+i(1-x)F_{m-1,i}(x)+(i-1)F_{m-1,i-1}(x) & \text{for\;} 1\leqslant i\leqslant m.
\end{cases}
\label{Eq:Recursive-fmi}
\end{align}
\end{thm}

In particular, for any integers $m,i$ with $1\leqslant i\leqslant m+1$, we have $F_{mi}(x)\in\mathbb{Z}[x]$.
From \eqref{Eq:Recursive-fmi}, we deduce that for any nonnegative integer $m$,
$$F_{m0}(x)=\frac{x}{2},\quad F_{m,m+1}(x)=-m!.$$
In general, we have the following result.

\begin{pro}
For any integers $m,i$ with $1\leqslant i\leqslant m+1$, we have $\deg F_{mi}(x)=m+1-i$, and the leading coefficient $c_{mi}$ of $F_{mi}(x)$ satisfies  $(-1)^{m+i}c_{mi}>0$.
\end{pro}

\pf We prove this result by induction on $m$. Assume that $m\geqslant 1$. The result for $i=m+1$ follows from $F_{m,m+1}(x)=-m!$. Now assume $1\leqslant i\leqslant m$, and
$$F_{m-1,i}(x)=c_{m-1,i}x^{m-i}+\text{lower degree terms}$$
with $(-1)^{m-1+i}c_{m-1,i}>0$.
Let $c_{m0}=\frac{1}{2}$. Then we obtain
$$F_{mi}(x)=(-ic_{m-1,i}+(i-1)c_{m-1,i-1})x^{m+1-i}+\text{lower degree terms},$$
and
\begin{align*}
&(-1)^{m+i}(-ic_{m-1,i}+(i-1)c_{m-1,i-1})\\
=&i(-1)^{m-1+i}c_{m-1,i}+(i-1)(-1)^{m-1+i-1}c_{m-1,i-1}>0,
\end{align*}
from which one can deduce the desired result. \hfill$\square$

Now we have
$$F_{m0}(x)=c_{m0}x$$
with $c_{m0}=\frac{1}{2}$, and for any integers $m,i$ with the condition $1\leqslant i\leqslant m+1$, we have
$$F_{mi}(x)=c_{mi}x^{m+1-i}+\text{lower degree terms},$$
with the recurrence relation
$$c_{mi}=-ic_{m-1,i}+(i-1)c_{m-1,i-1},\quad (1\leqslant i\leqslant m)$$
and $c_{m,m+1}=-m!$. In particular, according to the recurrence relation above, we deduce that if $m$ is a nonnegative integer, then $c_{m1}=(-1)^{m+1}$.

The following result will be used later.

\begin{pro}
For any nonnegative integer $m$, we have
\begin{align}
\sum\limits_{i=1}^{m+1}F_{mi}(x)x^{i-1}=-1,
\label{Eq:Sum-fmi}
\end{align}
and
\begin{align}
\sum\limits_{i=1}^{m+1}c_{mi}=-\delta_{m,0},
\label{Eq:Sum-cmi}
\end{align}
where $\delta_{i,j}$ is Kronecker's delta.
\end{pro}

\pf
We prove \eqref{Eq:Sum-fmi} by induction on $m$. The case of $m=0$ follows from the fact $F_{01}(t)=-1$. Now assume that $m\geqslant 1$, using the recurrence formula \eqref{Eq:Recursive-fmi}, we arrive at
\begin{align*}
\sum\limits_{i=1}^{m+1}F_{mi}(x)x^{i-1}=&\sum\limits_{i=1}^mF_{m-1,i}'(x)x^{i}+\sum\limits_{i=1}^{m}iF_{m-1,i}(x)x^{i-1}\\
&-\sum\limits_{i=1}^miF_{m-1,i}(x)x^i+\sum\limits_{i=1}^m(i-1)F_{m-1,i-1}(x)x^{i-1}-m!x^m\\
=&\sum\limits_{i=1}^m(F_{m-1,i}(x)x^{i})'-mF_{m-1,m}(x)x^m-m!x^m\\
=&\frac{d}{dx}\sum\limits_{i=1}^mF_{m-1,i}(x)x^{i}.
\end{align*}
Then we get \eqref{Eq:Sum-fmi} from the inductive hypothesis. Thus, comparing the coefficients of $x^m$ of both sides of \eqref{Eq:Sum-fmi}, we obtain the desired result \eqref{Eq:Sum-cmi}. \hfill$\square$

Now we use matrix computations to express $H(x)^i$ by $D^m F(x)$. First, for any nonnegative integer $m$, we define a $(m+1)\times (m+1)$ matrix $A_m(x)$ by
$$A_m(x)=\begin{pmatrix}
F_{01}(x) & &&\\
F_{11}(x) & F_{12}(x) &&\\
\vdots & \vdots & \ddots &\\
F_{m1}(x) & F_{m2}(x) & \cdots & F_{m,m+1}(x)
\end{pmatrix}.$$
It is clear that for $m\geqslant 1$, we have
$$A_m(x)=\begin{pmatrix}
A_{m-1}(x) & 0\\
\alpha_m(x) & -m!
\end{pmatrix}$$
with $\alpha_m(x)=(F_{m1}(x),\ldots,F_{mm}(x))$. Hence, the identity \eqref{Eq:Diff-f} can be rewritten as
\begin{align}
\begin{pmatrix}
F(x)\\
DF(x)\\
\vdots\\
D^mF(x)
\end{pmatrix}-\frac{1}{2}x\begin{pmatrix}
1\\
1\\
\vdots\\
1
\end{pmatrix}=A_m(x)\begin{pmatrix}
H(x)\\
H(x)^2\\
\vdots\\
H(x)^{m+1}
\end{pmatrix},
\label{Eq:G-h-Matrix}
\end{align}
and the identity \eqref{Eq:Sum-fmi} can be rewritten as
$$A_{m}(x)\begin{pmatrix}
1\\
x\\
x^2\\
\vdots\\
x^m
\end{pmatrix}=-\begin{pmatrix}
1\\
1\\
\vdots\\
1
\end{pmatrix}.$$
Since the matrix $A_m(x)$ is invertible, we find
\begin{align}
\begin{pmatrix}
H(x)\\
H(x)^2\\
\vdots\\
H(x)^{m+1}
\end{pmatrix}=A_m(x)^{-1}\begin{pmatrix}
F(x)\\
DF(x)\\
\vdots\\
D^mF(x)
\end{pmatrix}+\frac{1}{2}x\begin{pmatrix}
1\\
x\\
x^2\\
\vdots\\
x^m
\end{pmatrix}.\label{2.6}
\end{align}
Therefore, we need to obtain a description of $A_m(x)^{-1}$. From linear algebra, we know that the matrix $\begin{pmatrix}
A & 0\\
C & B
\end{pmatrix}$ is invertible with
$$\begin{pmatrix}
A & 0\\
C & B
\end{pmatrix}^{-1}=\begin{pmatrix}
A^{-1} & 0\\
-B^{-1}CA^{-1} & B^{-1}
\end{pmatrix},$$
provided that $A$ and $B$ are invertible square matrices. Hence by induction on $m$, we find that the inverses $A_m(x)^{-1}$ satisfy the recursive formula
\begin{align}
A_m(x)^{-1}=\begin{pmatrix}
A_{m-1}(x)^{-1} & 0\\
\frac{1}{m!}\alpha_m(x)A_{m-1}(x)^{-1} & \frac{-1}{m!}
\end{pmatrix},\quad (m\geqslant 1).
\label{Eq:Recursive-AmInverse}
\end{align}

For any nonnegative integer $m$, set
$$A_m(x)^{-1}=\begin{pmatrix}
G_{01}(x) & &&\\
G_{11}(x) & G_{12}(x) &&\\
\vdots & \vdots & \ddots &\\
G_{m1}(x) & G_{m2}(x) & \cdots & G_{m,m+1}(x)
\end{pmatrix}.$$
Then from (\ref{2.6}), for any positive integer $i$, we get
\begin{align}
H(x)^i=\sum\limits_{j=1}^iG_{i-1,j}(x)D^{j-1}F(x)+\frac{1}{2}x^i.
\label{Eq:h-i}
\end{align}
Moreover, it is easy to prove that for a nonnegative integer $m$, the functions $$1,F(x),D F(x),\ldots,D^m F(x)$$ are linearly independent over the rational function field $\mathbb{Q}(x)$. For a proof one can refer to \cite[Lemma 2.6]{LQ2019}.

Next, we give some properties of the polynomials  $G_{ij}(x)$.

\begin{pro} Let $m$ and $i$ be integers.
\begin{itemize}
  \item [(1)] For any $m\geqslant 0$, we have $G_{m,m+1}(x)=-\frac{1}{m!}$;
  \item [(2)] For $1\leqslant i\leqslant m$, we have the recursive formula
  \begin{align}
  G_{mi}(x)=\frac{1}{m!}\sum\limits_{j=i}^{m}F_{mj}(x)G_{j-1,i}(x);
  \label{Eq:Recursive-gmi}
  \end{align}
  \item [(3)] For $1\leqslant i\leqslant m+1$, we have $G_{mi}(x)\in\mathbb{Q}[x]$ with $\deg G_{mi}(x)\leqslant m+1-i$;
  \item [(4)] For $1\leqslant i\leqslant m+1$, set
  $$G_{mi}(x)=d_{mi}x^{m+1-i}+\text{lower degree terms}.$$
  Then we have $d_{m,m+1}=-\frac{1}{m!}$ and
   \begin{align}
  d_{mi}=\frac{1}{m!}\sum\limits_{j=i}^{m}c_{mj}d_{j-1,i}
  \label{Eq:Recursive-dmi}
  \end{align}
  for $1\leqslant i\leqslant m$.
\end{itemize}
\end{pro}

\pf The assertions in items (1) and (2) follow from \eqref{Eq:Recursive-AmInverse}. To prove the item (3), we proceed by induction on $m$. For the case of $m=0$, we get the result from $G_{01}(t)=-1$. Assume that $m\geqslant 1$, then $G_{m,m+1}(x)=-\frac{1}{m!}\in\mathbb{Q}[x]$ with degree zero. For $1\leqslant i\leqslant j\leqslant m$, by the induction assumption, we may set
$$G_{j-1,i}(x)=d_{j-1,i}x^{j-i}+\text{lower degree terms}\in\mathbb{Q}[x].$$
Since
$$F_{mj}(x)=c_{mj}x^{m+1-j}+\text{lower degree terms}\in\mathbb{Z}[x],$$
we get
$$F_{mj}(x)G_{j-1,i}(x)=c_{mj}d_{j-1,i}x^{m+1-i}+\text{lower degree terms}\in\mathbb{Q}[x].$$
Using \eqref{Eq:Recursive-gmi}, we finally get
$$G_{mi}(x)=\left(\frac{1}{m!}\sum\limits_{j=i}^{m}c_{mj}d_{j-1,i}\right)x^{m+1-i}+\text{lower degree terms}\in\mathbb{Q}[x].$$
The item (4) follows from the above proof.\hfill$\square$
\begin{cor}
For any nonnegative integer $m$, we have $d_{m1}=-1$.
\end{cor}
\pf
We use induction on $m$. If $m\geqslant 1$, using \eqref{Eq:Recursive-dmi} and the induction assumption, we get
$$d_{m1}=-\frac{1}{m!}\sum\limits_{j=1}^mc_{mj}.$$
By \eqref{Eq:Sum-cmi}, we have
$$d_{m1}=-\frac{1}{m!}(-\delta_{m,0}-c_{m,m+1}),$$
which implies the result.\hfill$\square$

\subsection{A weighted sum formula of the Bernoulli numbers}

Let $n$ be a fixed positive integer, and $m_1,\ldots,m_n$ be fixed nonnegative integers. Set $|{\mathbf{m}}|_n:=m_1+m_2+\cdots+m_n+n$. Now, we evaluate $D^{m_1}F(x)\cdots D^{m_n}F(x)$. First, using the fact that
$$D^mF(x)=\sum\limits_{i=1}^{\infty}(2i)^m\frac{\beta_{2i}}{(2i)!}x^{2i}=\sum\limits_{i=0}^{\infty}(2i)^m\frac{\beta_{2i}}{(2i)!}x^{2i},$$
we get
\begin{align*}
D^{m_1}F(x)\cdots D^{m_n}F(x)=\sum\limits_{k=n}^\infty \left( \sum\limits_{(k_1,\ldots,k_n)\in I(k,n)} (2k_1)^{m_1}\cdots (2k_n)^{m_n} \frac{\beta_{2k_1}\cdots \beta_{2k_n}}{(2k_1)!\cdots (2k_n)!}\right)x^{2k}.
\end{align*}
Therefore, for any positive integer $k$ with $k\geqslant n$, the coefficient of $x^{2k}$ in $D^{m_1}F(x)\cdots D^{m_n}F(x)$ is
\begin{align}
\sum\limits_{(k_1,\ldots,k_n)\in I(k,n)}(2k_1)^{m_1}\cdots(2k_n)^{m_n}\frac{\beta_{2k_1}\cdots \beta_{2k_n}}{(2k_1)!\cdots(2k_n)!}.
\label{Eq:Coeff-Left}
\end{align}

Next, using \eqref{Eq:Diff-f}, we have
$$D^{m_1}F(x)\cdots D^{m_n}F(x)=\sum\limits_{i=0}^{|\mathbf{m}|_n}F_i(x)H(x)^i,$$
with
$$F_i(x)=\sum\limits_{i_1+\cdots+i_n=i\atop 0\leqslant i_j\leqslant m_j+1}F_{m_1i_1}(x)\cdots F_{m_ni_n}(x).$$

\begin{pro}
We have
$$F_0(x)=\left(\frac{x}{2}\right)^n,$$
and $\deg F_i(x)\leqslant |\mathbf{m}|_n-i$ for any nonnegative integer $i$.
\end{pro}

\pf
For integers $i_1,\ldots,i_n$ with the conditions $i_1+\cdots+i_n=i$ and $0\leqslant i_j\leqslant m_j+1$, we have
$$\deg(F_{m_1i_1}(x)\cdots F_{m_ni_n}(x))\leqslant \sum\limits_{j=1}^n(m_j+1-i_j)=|\mathbf{m}|_n-i,$$
which implies that $\deg F_i(x)\leqslant|\mathbf{m}|_n-i$.\hfill$\square$

Then using \eqref{Eq:h-i}, we get
\begin{align*}
&D^{m_1}F(x)\cdots D^{m_n}F(x)\\
=&\sum\limits_{i=1}^{|\mathbf{m}|_n}F_i(x)\left(\sum\limits_{j=1}^iG_{i-1,j}(x)D^{j-1}F(x)+\frac{1}{2}x^i\right)+R_0(x)\\
=&\sum\limits_{j=1}^{|\mathbf{m}|_n}R_j(x)D^{j-1}F(x)+R_0(x)
\end{align*}
with
$$R_0(x):=F_0(x)+\frac{1}{2}\sum\limits_{i=1}^{|\mathbf{m}|_n}F_i(x)x^i$$
and
$$R_j(x):=\sum\limits_{i=j}^{|\mathbf{m}|_n}F_i(x)G_{i-1,j}(x),\quad (1\leqslant j\leqslant |\mathbf{m}|_n).$$

\begin{pro}\label{pro7}
Let $j$ be a nonnegative integer with $j\leqslant |\mathbf{m}|_n$. Then
\begin{itemize}
  \item [(1)] the function $R_j(x)$ is even;
  \item [(2)] we have
  $$R_0(x)=\frac{1}{2^{n+1}}\left(x^n+(-x)^n \right).$$
In particular, $\deg R_0(x)\leqslant n$;
  \item [(3)] for $j>0$, we have $\deg R_j(x)\leqslant |\mathbf{m}|_n-j$.
Moreover, we have $\deg R_1(x)\leqslant |\mathbf{m}|_n-2$ provided that $n$ is even or $m_1,\ldots,m_n$ are not all zero.
\end{itemize}
\end{pro}

\pf
Since $D^{m}F(x)$ is even, we have
$$\sum\limits_{j=1}^{|\mathbf{m}|_n}R_j(x)D^{j-1}F(x)+R_0(x)=\sum\limits_{j=1}^{|\mathbf{m}|_n}R_j(-x)D^{j-1}F(x)+R_0(-x).$$
Using the fact that the functions $1,F(x),D F(x),\ldots,D^m F(x)$ are linearly independent over the rational function field $\mathbb{Q}(x)$, we know all $R_j(x)$ are even functions.

By the definition of $F_i(x)$, we have
$$\sum\limits_{i=0}^{|\mathbf{m}|_n}F_i(x)x^i=\prod\limits_{j=1}^n\sum\limits_{i_j=0}^{m_j+1}F_{m_ji_j}(x)x^{i_j}.$$
Using \eqref{Eq:Sum-fmi}, we find
$$\sum\limits_{i=0}^{|\mathbf{m}|_n}F_i(x)x^i=\prod\limits_{j=1}^n(F_{m_j0}(x)-x).$$
Then we get (2) from the fact that $F_{m0}(x)=\frac{1}{2}x$ and the expression of $F_0(x)$.

Since
$$\deg F_i(x)G_{i-1,j}(x)\leqslant (|\mathbf{m}|_n-i)+(i-j)=|\mathbf{m}|_n-j,$$
we get $\deg R_j(x)\leqslant |\mathbf{m}|-j$.

If we set
$$\widetilde{c}_{mi}=\begin{cases}
\frac{1}{2}\delta_{m,0} & \text{if\;} i=0,\\
c_{mi} & \text{if\;} i\neq 0,
\end{cases}$$
then the coefficient of $x^{m+1-i}$ in $F_{mi}(x)$ is $\widetilde{c}_{mi}$ for any integers $m,i$ with the condition $0\leqslant i\leqslant m+1$.
Since
$$R_1(x)=\sum\limits_{i=1}^{|\mathbf{m}|_n}\sum\limits_{i_1+\cdots+i_n=i\atop 0\leqslant i_j\leqslant m_j+1}F_{m_1i_1}(x)\cdots F_{m_ni_n}(t)G_{i-1,1}(x),$$
and $d_{i-1,1}=-1$, we find the coefficient of $x^{|\mathbf{m}|-1}$ in $R_1(x)$ is
\begin{align*}
&-\sum\limits_{i=1}^{|\mathbf{m}|_n}\sum\limits_{i_1+\cdots+i_n=i\atop 0\leqslant i_j\leqslant m_j+1}\widetilde{c}_{m_1i_1}\cdots \widetilde{c}_{m_ni_n}
=\widetilde{c}_{m_10}\cdots \widetilde{c}_{m_n0}-\prod\limits_{j=1}^n\sum\limits_{i_j=0}^{m_j+1}\widetilde{c}_{m_ji_j},
\end{align*}
more precisely which equals
$$\widetilde{c}_{m_10}\cdots \widetilde{c}_{m_n0}-\prod\limits_{j=1}^n(\widetilde{c}_{m_j0}-\delta_{m_j,0})$$
by \eqref{Eq:Sum-cmi}. Then the coefficient of $x^{|\mathbf{m}|_n-1}$ in $R_1(x)$ is
$$\left(\frac{1}{2}\right)^n(1-(-1)^n)\delta_{m_1,0}\cdots\delta_{m_n,0},$$
which is zero if $n$ is even or at least one $m_i$ is not zero. \hfill$\square$

Let $a_{jl}\in \mathbb{Q}$ be the coefficient of $x^{2l}$ in the even polynomial $R_j(x)$, then we have
\begin{align}
R_j(x)=\sum\limits_{l\geqslant 0}a_{jl}x^{2l}=\sum\limits_{l=0}^{\left[(|\mathbf{m}|_n-j)/{2}\right]}a_{jl}x^{2l}.
\label{Eq:Fj}
\end{align}
Moreover, from Proposition \ref{pro7}, it is clear that if
$n$ is even or $m_1,\ldots,m_n$ are not all zero, then
$$R_1(x)=\sum\limits_{l=0}^{\left[(|\mathbf{m}|_n-2)/{2}\right]}a_{1l}x^{2l}.$$
Hence we have
$$D^{m_1}F(x)\cdots D^{m_n}F(x)=\sum\limits_{j=1}^{|\mathbf{m}|_n}\sum\limits_{l=0}^{\left[(|\mathbf{m}|_n-j)/{2}\right]}a_{jl}x^{2l}D^{j-1}F(x)+R_0(x).$$
Changing the order of the summation yields
$$D^{m_1}F(x)\cdots D^{m_n}F(x)=\sum\limits_{l=0}^{T}\sum\limits_{j=1}^{|\mathbf{m}|_n-2l}a_{jl}x^{2l}D^{j-1}F(x)+R_0(x),$$
where
$$T=\begin{cases}
\left[\frac{n-1}{2}\right] & \text{if\;} m_1=\cdots=m_n=0,\\
&\\
\left[\frac{|\mathbf{m}|_n-2}{2}\right] & \text{otherwise}.
\end{cases}$$
Since
$$D^{j-1}F(x)=\sum\limits_{i=0}^\infty(2i)^{j-1}\frac{\beta_{2i}}{(2i)!}x^{2i},$$
we get
\begin{align*}
&D^{m_1}F(x)\cdots D^{m_n}F(x)
=&\sum\limits_{k=0}^\infty\sum\limits_{l=0}^{\min\{T,k\}}\left(\sum\limits_{j=1}^{|\mathbf{m}|_n-2l}a_{jl}(2k-2l)^{j-1}\right)\frac{\beta_{2k-2l}}{(2k-2l)!}x^{2k}+R_0(x).
\end{align*}
Then the coefficient of $x^{2k}$ in $D^{m_1}F(x)\cdots D^{m_n}F(x)$ is
\begin{align}
\sum\limits_{l=0}^{\min\{T,k\}}\left(\sum\limits_{j=1}^{|\mathbf{m}|_n-2l}2^{j-1}a_{jl}(k-l)^{j-1}\right)\frac{\beta_{2k-2l}}{(2k-2l)!},
\label{Eq:Coeff-Right}
\end{align}
provided that $k\geqslant n$.

Finally, comparing \eqref{Eq:Coeff-Right} with \eqref{Eq:Coeff-Left}, we get a weighted sum formula of the Bernoulli numbers.

\begin{thm}\label{Thm:WeightedSum-Bernoulli}
Let $n,k$ be positive integers with $k\geqslant n$. Then for any nonnegative integers $m_1,\ldots,m_n$, we have
\begin{align}
&\sum\limits_{k_1+\cdots+k_n=k\atop k_j\geqslant 1}k_1^{m_1}\cdots k_n^{m_n}\frac{\beta_{2k_1}\cdots \beta_{2k_n}}{(2k_1)!\cdots(2k_n)!}
=&\sum\limits_{l=0}^{\min\{T,k\}}\left(\sum\limits_{j=1}^{|\mathbf{m}|_n-2l}\frac{a_{jl}(k-l)^{j-1}}{2^{m_1+\cdots+m_n-j+1}}\right)\frac{\beta_{2k-2l}}{(2k-2l)!},
\label{Eq:WeightedSum-Bernoulli}
\end{align}
where $T=\max\{[(|\mathbf{m}|_n-2)/2],[(n-1)/2]\}$ and $a_{jl}$ are determined by \eqref{Eq:Fj}.
\end{thm}

\subsection{Proof of Theorem \ref{thm1}}

Now, we prove the weighted sum formula \eqref{1.1} of $t$-values at even arguments. Using Euler's formula of $\zeta(2k)$, we have
\begin{align}
t(2k)=(-1)^{k+1}\frac{\beta_{2k}}{2(2k)!}\pi^{2k}.
\label{Eq:Euler-Formula}
\end{align}
Then from Theorem \ref{Thm:WeightedSum-Bernoulli}, we get the following weighted sum formula of $t$-values at even arguments.

\begin{thm}\label{Thm:WeightedSum-Zeta}
Let $n,k$ be positive integers with $k\geqslant n$. Then for any nonnegative integers $m_1,\ldots,m_n$, we have
\begin{align}
&\sum\limits_{k_1+\cdots+k_n=k\atop k_j\geqslant 1}k_1^{m_1}\cdots k_n^{m_n}t(2k_1)\cdots t(2k_n)\nonumber\\ &=(-1)^n\sum\limits_{l=0}^{\min\{T,k\}}\frac{(2l)!}{B_{2l}}\left(\sum\limits_{j=1}^{|\mathbf{m}|_n-2l}\frac{a_{jl}(k-l)^{j-1}}{2^{|\mathbf{m}|_n+2l-j-1}}\right)\zeta(2l)t(2k-2l),
\label{Eq:WeightedSum-Zeta}
\end{align}
where $\zeta(0)=-1/2$ and $t(0)=0$, $T=\max\{[(|\mathbf{m}|_n-2)/2],[(n-1)/2]\}$ and $a_{jl}$ are determined by \eqref{Eq:Fj}.
\end{thm}

Finally, from Theorem \ref{Thm:WeightedSum-Zeta}, we prove Theorem \ref{thm1}.\hfill$\square$

\subsection{Proof of Theorem \ref{thm2}}

Next, we use the symmetric sum formulas of Hoffman \cite[Theorems 2.5 and 2.8]{H2016} to prove Theorem \ref{thm2}. For a partition $\Pi=\{P_1,P_2,\ldots,P_i\}$ of the set $\{1,2,\ldots,n\}$, let $l_j=\sharp P_j$ and
$$c(\Pi)=\prod\limits_{j=1}^i (l_j-1)!,\quad \tilde{c}(\Pi)=(-1)^{n-i}c(\Pi).$$
We also denote by $\mathcal{P}_n$ the set of all partitions of the set $\{1,2,\ldots,n\}$. Then the symmetric sum formulas of multiple $t$-values are
\begin{align}
\sum\limits_{\sigma\in S_n}t(k_{\sigma(1)},\ldots,k_{\sigma(n)})=\sum\limits_{\Pi\in\mathcal{P}_n}\tilde{c}(\Pi)t(\mathbf{k},\Pi)
\label{Eq:SymSum-MZV}
\end{align}
and
\begin{align}
\sum\limits_{\sigma\in S_n}t^{\star}(k_{\sigma(1)},\ldots,k_{\sigma(n)})=\sum\limits_{\Pi\in\mathcal{P}_n}c(\Pi)t(\mathbf{k},\Pi),
\label{Eq:SymSum-MZSV}
\end{align}
where $\mathbf{k}=(k_1,\ldots,k_n)$ is an index with all $k_i>1$, $S_n$ is the symmetric group of degree $n$ and for a partition $\Pi=\{P_1,\ldots,P_i\}\in\mathcal{P}_n$,
$$t(\mathbf{k},\Pi)=\prod\limits_{j=1}^i t\left(\sum\limits_{l\in P_j}k_l\right).$$
Now let $\mathbf{k}=(2k_1,\ldots,2k_n)$ with all $k_i$ positive integers. Using \eqref{Eq:SymSum-MZV} and \eqref{Eq:SymSum-MZSV}, we have
\begin{align}
&\sum\limits_{\sigma\in S_n}t(2k_{\sigma(1)},\ldots,2k_{\sigma(n)})\nonumber\\
=&\sum\limits_{i=1}^n(-1)^{n-i}\sum\limits_{l_1+\cdots+l_i=n\atop l_1\geqslant \cdots\geqslant l_i\geqslant 1}\prod\limits_{j=1}^i(l_j-1)!\sum\limits_{\Pi=\{P_1,\ldots,P_i\}\in\mathcal{P}_n\atop \sharp{P_j}=l_j}t(\mathbf{k},\Pi)
\label{Eq:SymSum-2-MZV}
\end{align}
and
\begin{align}
&\sum\limits_{\sigma\in S_n}t^{\star}(2k_{\sigma(1)},\ldots,2k_{\sigma(n)})\nonumber\\
=&\sum\limits_{i=1}^n\sum\limits_{l_1+\cdots+l_i=n\atop l_1\geqslant \cdots\geqslant l_i\geqslant 1 }\prod\limits_{j=1}^i(l_j-1)!\sum\limits_{\Pi=\{P_1,\ldots,P_i\}\in\mathcal{P}_n\atop \sharp{P_j}=l_j}t(\mathbf{k},\Pi).
\label{Eq:SymSum-2-MZSV}
\end{align}

From now on, let $k,n$ be fixed positive integers with $k\geqslant n$, and let $f(x_1,\ldots,x_n)$ be a fixed symmetric polynomial with rational coefficients. It is easy to see that
\begin{align*}
&\sum\limits_{(k_1,\ldots,k_n)\in I(k,n)}f(k_1,\ldots,k_n)\sum\limits_{\sigma\in S_n}t(2k_{\sigma(1)},\ldots,2k_{\sigma(n)})\\
=&n!\sum\limits_{(k_1,\ldots,k_n)\in I(k,n)}f(k_1,\ldots,k_n)t(2k_1,\ldots,2k_n)=n!T_f(2,k,n)
\end{align*}
and
\begin{align*}
&\sum\limits_{(k_1,\ldots,k_n)\in I(k,n)}f(k_1,\ldots,k_n)\sum\limits_{\sigma\in S_n}t^{\star}(2k_{\sigma(1)},\ldots,2k_{\sigma(n)})\\
=&n!\sum\limits_{(k_1,\ldots,k_n)\in I(k,n)}f(k_1,\ldots,k_n)t^{\star}(2k_1,\ldots,2k_n)=n!T_f^\star(2,k,n).
\end{align*}
On the other hand, for a partition $\Pi=\{P_1,\ldots,P_i\}\in\mathcal{P}_n$ with $\sharp P_j=l_j$, we have
\begin{align}
&\sum\limits_{(k_1,\ldots,k_n)\in I(k,n)}f(k_1,\ldots,k_n)t(\mathbf{k},\Pi)\nonumber\\
=&\sum\limits_{s_1+\cdots+s_i=k\atop s_j\geqslant 1}\sum\limits_{{{k_1+\cdots+k_{l_1}=s_1\atop\vdots}\atop k_{l_1+\cdots+l_{i-1}+1}+\cdots+k_n=s_i}\atop k_j\geqslant 1}f(k_1,\ldots,k_n)t(2s_1)\cdots t(2s_i).
\label{Eq:F-times-zeta}
\end{align}
To treat the inner sum about $f(k_1,\ldots,k_n)$ in the right-hand side of \eqref{Eq:F-times-zeta}, we need the following lemma.

\begin{lem}[{\cite[Lemma 4.2]{LQ2019}}]\label{Lem:PowerSum}
Let $k$ and $n$ be integers with $k\geqslant n\geqslant 1$, and let $p_1,\ldots,p_n$ be nonnegative integers.  Then there exists a polynomial $g(x)\in\mathbb{Q}[x]$ of degree $p_1+\cdots+p_n+n-1$, such that
$$\sum\limits_{k_1+\cdots+k_n=k\atop k_j\geqslant 1}k_1^{p_1}\cdots k_n^{p_n}=g(k).$$
\end{lem}

 Using Lemma \ref{Lem:PowerSum}, there exists a polynomial $g_{s_1,\ldots,s_i}(x_1,\ldots,x_i)\in\mathbb{Q}[x_1,\ldots,x_i]$ of degree $\deg f+n-i$, such that
\begin{align*}
&\sum\limits_{(k_1,\ldots,k_n)\in I(k,n)}f(k_1,\ldots,k_n)t(\mathbf{k},\Pi)\\
=&\sum\limits_{s_1+\cdots+s_i=k\atop s_j\geqslant 1}g_{s_1,\ldots,s_i}(s_1,\ldots,s_i)t(2s_1)\cdots t(2s_i).
\end{align*}
Therefore  we get
\begin{align*}
T_f(2,k,n)=&\frac{1}{n!}\sum\limits_{i=1}^n(-1)^{n-i}\sum\limits_{l_1+\cdots+l_i=n\atop l_1\geqslant\cdots \geqslant l_i\geqslant 1}\prod\limits_{j=1}^i(l_j-1)!n(l_1,\ldots,l_i)\\
&\times\sum\limits_{s_1+\cdots+s_i=k\atop s_j\geqslant 1}g_{s_1,\ldots,s_i}(s_1,\ldots,s_i)t(2s_1)\cdots t(2s_i)
\end{align*}
and
\begin{align*}
T^{\star}_f(2,k,n)=&\frac{1}{n!}\sum\limits_{i=1}^n\sum\limits_{l_1+\cdots+l_i=n\atop l_1\geqslant\cdots\geqslant l_i\geqslant 1}\prod\limits_{j=1}^i(l_j-1)!n(l_1,\ldots,l_i)\\
&\times\sum\limits_{s_1+\cdots+s_i=k\atop s_j\geqslant 1}g_{s_1,\ldots,s_i}(s_1,\ldots,s_i)t(2s_1)\cdots t(2s_i),
\end{align*}
where
$$n(l_1,\ldots,l_i)=\frac{n!}{\prod\limits_{j=1}^il_j!\prod\limits_{j=1}^n\sharp\{m\mid 1\leqslant m\leqslant i,k_m=j\}!}$$
is the number of partitions $\Pi=\{P_1,\ldots,P_i\}\in\mathcal{P}_n$ with the conditions $\sharp P_j=l_j$ for $j=1,2,\ldots,i$.

Thus, applying Theorem \ref{thm1}, we prove the weighted sum formulas (\ref{1.2}) and (\ref{1.3}). \hfill$\square$

%%%%%%%%%%%%%%%%%%%%%%%%%%%%%%%%%%%%%%%%%%%%%%%%%%%%%%%%%%%%%%%%%%%%%%%%%%%%%%%%%%%%%%%%%%%%%%%%%%%%%%%%%%

\appendix

\section{Some weighted sum formulas through depth $4$}\label{Sec:Example}

In this appendix, we list some explicit weighted sum formulas of depth $n\leqslant 4$. For any positive integers $k,n$ with $k\geqslant n$, we set
$$\sum\nolimits^{(n)}=\sum_{(k_1,\ldots,k_n)\in I(k,n)}.$$

\subsection{Weighted sum formulas of the Bernoulli numbers}

Recall that $\beta_{2i}=(2^{2i}-1)B_{2i}$. If $n=2$, we have
\begin{align*}
&\sum\nolimits^{(2)} \frac{\beta_{2k_1}\beta_{2k_2}}{(2k_1)!(2k_2)!}=-(2k-1)\frac{\beta_{2k}}{(2k)!},\\
&\sum\nolimits^{(2)} k_1\frac{\beta_{2k_1}\beta_{2k_2}}{(2k_1)!(2k_2)!}=-\frac{1}{2}k(2k-1)\frac{\beta_{2k}}{(2k)!},\\
&\sum\nolimits^{(2)} k_1^2\frac{\beta_{2k_1}\beta_{2k_2}}{(2k_1)!(2k_2)!}=-\frac{1}{12}k(2k-1)(4k-1)\frac{\beta_{2k}}{(2k)!}-\frac{1}{24}(2k-3)\frac{\beta_{2k-2}}{(2k-2)!},\\
&\sum\nolimits^{(2)} k_1 k_2\frac{\beta_{2k_1}\beta_{2k_2}}{(2k_1)!(2k_2)!}=-\frac{1}{12}k(2k-1)(2k+1)\frac{\beta_{2k}}{(2k)!}+\frac{1}{24}(2k-3)\frac{\beta_{2k-2}}{(2k-2)!},\\
&\sum\nolimits^{(2)} k_1^3\frac{\beta_{2k_1}\beta_{2k_2}}{(2k_1)!(2k_2)!}=-\frac{1}{8}k^2(2k-1)^2\frac{\beta_{2k}}{(2k)!}-\frac{1}{16}k(2k-3)\frac{\beta_{2k-2}}{(2k-2)!},\\
&\sum\nolimits^{(2)} k_1^2k_2\frac{\beta_{2k_1}\beta_{2k_2}}{(2k_1)!(2k_2)!}=-\frac{1}{24}k^2(2k-1)(2k+1)\frac{\beta_{2k}}{(2k)!}+\frac{1}{48}k(2k-3)\frac{\beta_{2k-2}}{(2k-2)!},\\
&\sum\nolimits^{(2)} k_1^4\frac{\beta_{2k_1}\beta_{2k_2}}{(2k_1)!(2k_2)!}=-\frac{1}{240}k(2k-1)(4k-1)(12k^2-6k-1)\frac{\beta_{2k}}{(2k)!}\\
&\qquad\quad -\frac{1}{96}(2k-3)(8k^2-6k+5)\frac{\beta_{2k-2}}{(2k-2)!}+\frac{1}{480}(2k-5)\frac{\beta_{2k-4}}{(2k-4)!},\\
&\sum\nolimits^{(2)} k_1^3k_2\frac{\beta_{2k_1}\beta_{2k_2}}{(2k_1)!(2k_2)!}=-\frac{1}{240}k(2k-1)(2k+1)(6k^2-1)\frac{\beta_{2k}}{(2k)!}\\
&\qquad\quad +\frac{1}{96}(2k-3)(2k^2-6k+5)\frac{\beta_{2k-2}}{(2k-2)!}-\frac{1}{480}(2k-5)\frac{\beta_{2k-4}}{(2k-4)!},\\
&\sum\nolimits^{(2)} k_1^2k_2^2\frac{\beta_{2k_1}\beta_{2k_2}}{(2k_1)!(2k_2)!}=-\frac{1}{240}k(2k-1)(2k+1)(4k^2+1)\frac{\beta_{2k}}{(2k)!}\\
&\qquad\quad +\frac{1}{96}(2k-3)(6k-5)\frac{\beta_{2k-2}}{(2k-2)!}+\frac{1}{480}(2k-5)\frac{\beta_{2k-4}}{(2k-4)!},\\
&\sum\nolimits^{(2)} k_1^5\frac{\beta_{2k_1}\beta_{2k_2}}{(2k_1)!(2k_2)!}=-\frac{1}{96}k^2(2k-1)^2(8k^2-4k-1)\frac{\beta_{2k}}{(2k)!}\\
&\qquad\quad -\frac{5}{192}k(2k-3)(4k^2-6k+5)\frac{\beta_{2k-2}}{(2k-2)!}+\frac{1}{192}k(2k-5)\frac{\beta_{2k-4}}{(2k-4)!},\\
&\sum\nolimits^{(2)} k_1^4k_2\frac{\beta_{2k_1}\beta_{2k_2}}{(2k_1)!(2k_2)!}=-\frac{1}{480}k^2(2k-1)(2k+1)(8k^2-3)\frac{\beta_{2k}}{(2k)!}\\
&\qquad\quad +\frac{1}{192}k(2k-3)(4k^2-18k+15)\frac{\beta_{2k-2}}{(2k-2)!}-\frac{1}{320}k(2k-5)\frac{\beta_{2k-4}}{(2k-4)!},\\
&\sum\nolimits^{(2)} k_1^3k_2^2\frac{\beta_{2k_1}\beta_{2k_2}}{(2k_1)!(2k_2)!}=-\frac{1}{480}k^2(2k-1)(2k+1)(4k^2+1)\frac{\beta_{2k}}{(2k)!}\\
&\qquad\quad +\frac{1}{192}k(2k-3)(6k-5)\frac{\beta_{2k-2}}{(2k-2)!}+\frac{1}{960}k(2k-5)\frac{\beta_{2k-4}}{(2k-4)!}.
\end{align*}

If $n=3$, we have
\begin{align*}
&\sum\nolimits^{(3)}\frac{\beta_{2k_1}\beta_{2k_2}\beta_{2k_3}}{(2k_1)!(2k_2)!(2k_3)!}=(k-1)(2k-1)\frac{\beta_{2k}}{(2k)!}+\frac{1}{4}\frac{\beta_{2k-2}}{(2k-2)!},\\
&\sum\nolimits^{(3)} k_1\frac{\beta_{2k_1}\beta_{2k_2}\beta_{2k_3}}{(2k_1)!(2k_2)!(2k_3)!}=\frac{1}{3}k(k-1)(2k-1)\frac{\beta_{2k}}{(2k)!}+\frac{1}{12}k\frac{\beta_{2k-2}}{(2k-2)!},\\
&\sum\nolimits^{(3)} k_1^2\frac{\beta_{2k_1}\beta_{2k_2}\beta_{2k_3}}{(2k_1)!(2k_2)!(2k_3)!}=\frac{1}{12}k(k-1)(2k-1)^2\frac{\beta_{2k}}{(2k)!}+\frac{1}{24}(4k^2-11k+9)\frac{\beta_{2k-2}}{(2k-2)!},\\
&\sum\nolimits^{(3)} k_1k_2\frac{\beta_{2k_1}\beta_{2k_2}\beta_{2k_3}}{(2k_1)!(2k_2)!(2k_3)!}=\frac{1}{24}k(k-1)(2k-1)(2k+1)\frac{\beta_{2k}}{(2k)!}\\
&\qquad\qquad\qquad\qquad\qquad\qquad -\frac{1}{48}(k-1)(2k-9)\frac{\beta_{2k-2}}{(2k-2)!},\\
&\sum\nolimits^{(3)} k_1^3\frac{\beta_{2k_1}\beta_{2k_2}\beta_{2k_3}}{(2k_1)!(2k_2)!(2k_3)!}=\frac{1}{120}k(k-1)(2k-1)(12k^2-12k+1)\frac{\beta_{2k}}{(2k)!}\\
&\qquad\qquad +\frac{1}{48}(8k^3-24k^2+17k+3)\frac{\beta_{2k-2}}{(2k-2)!}+\frac{1}{240}(2k-5)\frac{\beta_{2k-4}}{(2k-4)!},\\
&\sum\nolimits^{(3)} k_1^2k_2\frac{\beta_{2k_1}\beta_{2k_2}\beta_{2k_3}}{(2k_1)!(2k_2)!(2k_3)!}=\frac{1}{240}k(k-1)(2k-1)(2k+1)(4k-1)\frac{\beta_{2k}}{(2k)!}\\
&\qquad\qquad +\frac{1}{96}(k-1)(2k+3)\frac{\beta_{2k-2}}{(2k-2)!}-\frac{1}{480}(2k-5)\frac{\beta_{2k-4}}{(2k-4)!},\\
&\sum\nolimits^{(3)} k_1k_2k_3\frac{\beta_{2k_1}\beta_{2k_2}\beta_{2k_3}}{(2k_1)!(2k_2)!(2k_3)!}=\frac{1}{120}k(k-1)(k+1)(2k-1)(2k+1)\frac{\beta_{2k}}{(2k)!}\\
&\qquad\qquad -\frac{1}{48}(k-3)(k-1)(2k-1)\frac{\beta_{2k-2}}{(2k-2)!}+\frac{1}{240}(2k-5)\frac{\beta_{2k-4}}{(2k-4)!},\\
&\sum\nolimits^{(3)} k_1^4\frac{\beta_{2k_1}\beta_{2k_2}\beta_{2k_3}}{(2k_1)!(2k_2)!(2k_3)!}=\frac{1}{240}k(k-1)(2k-1)^2(8k^2-8k-1)\frac{\beta_{2k}}{(2k)!}\\
&\qquad\qquad+\frac{1}{96}(16k^4-64k^3+96k^2-79k+39)\frac{\beta_{2k-2}}{(2k-2)!}\\
&\qquad\qquad +\frac{1}{480}(2k-5)(3k+1)\frac{\beta_{2k-4}}{(2k-4)!},\\
&\sum\nolimits^{(3)} k_1^3k_2\frac{\beta_{2k_1}\beta_{2k_2}\beta_{2k_3}}{(2k_1)!(2k_2)!(2k_3)!}=\frac{1}{480}k(k-1)(2k-1)(2k+1)(4k^2-2k-1)\frac{\beta_{2k}}{(2k)!}\\
&\qquad+\frac{1}{192}(k-1)(16k^2-46k+39)\frac{\beta_{2k-2}}{(2k-2)!}-\frac{1}{960}(k+1)(2k-5)\frac{\beta_{2k-4}}{(2k-4)!},\\
&\sum\nolimits^{(3)} k_1^2k_2^2\frac{\beta_{2k_1}\beta_{2k_2}\beta_{2k_3}}{(2k_1)!(2k_2)!(2k_3)!}=\frac{1}{1440}k(k-1)(2k-1)(2k+1)(8k^2-4k+3)\frac{\beta_{2k}}{(2k)!}\\
&\qquad\qquad+\frac{1}{576}(k-1)(8k^3-64k^2+168k-117)\frac{\beta_{2k-2}}{(2k-2)!}\\
&\qquad\qquad-\frac{1}{2880}(2k-5)(7k-3)\frac{\beta_{2k-4}}{(2k-4)!},\\
&\sum\nolimits^{(3)} k_1^2k_2k_3\frac{\beta_{2k_1}\beta_{2k_2}\beta_{2k_3}}{(2k_1)!(2k_2)!(2k_3)!}=\frac{1}{360}k^2(k-1)(k+1)(2k-1)(2k+1)\frac{\beta_{2k}}{(2k)!}\\
&\qquad\quad -\frac{1}{144}k(k-3)(k-1)(2k-1)\frac{\beta_{2k-2}}{(2k-2)!}+\frac{1}{720}k(2k-5)\frac{\beta_{2k-4}}{(2k-4)!}.
\end{align*}

If $n=4$, we have
\begin{align*}
&\sum\nolimits^{(4)}\frac{\beta_{2k_1}\beta_{2k_2}\beta_{2k_3}\beta_{2k_4}}{(2k_1)!(2k_2)!(2k_3)!(2k_4)!}=-\frac{1}{3}(k-1)(2k-3)(2k-1)\frac{\beta_{2k}}{(2k)!}\\
&\qquad\qquad\qquad\qquad\qquad\qquad\qquad\qquad -\frac{1}{3}(2k-3)\frac{\beta_{2k-2}}{(2k-2)!},\\
&\sum\nolimits^{(4)} k_1\frac{\beta_{2k_1}\beta_{2k_2}\beta_{2k_3}\beta_{2k_4}}{(2k_1)!(2k_2)!(2k_3)!(2k_4)!}=-\frac{1}{12}k(k-1)(2k-3)(2k-1)\frac{\beta_{2k}}{(2k)!}\\
&\qquad\qquad\qquad\qquad\qquad\qquad\qquad -\frac{1}{12}k(2k-3)\frac{\beta_{2k-2}}{(2k-2)!},\\
&\sum\nolimits^{(4)} k_1^2\frac{\beta_{2k_1}\beta_{2k_2}\beta_{2k_3}\beta_{2k_4}}{(2k_1)!(2k_2)!(2k_3)!(2k_4)!}=-\frac{1}{120}k(k-1)(2k-3)(2k-1)(4k-3)\frac{\beta_{2k}}{(2k)!}\\
&\qquad\qquad -\frac{1}{12}(2k-3)(k^2-3k+3)\frac{\beta_{2k-2}}{(2k-2)!}-\frac{1}{160}(2k-5)\frac{\beta_{2k-4}}{(2k-4)!},\\
&\sum\nolimits^{(4)} k_1k_2\frac{\beta_{2k_1}\beta_{2k_2}\beta_{2k_3}\beta_{2k_4}}{(2k_1)!(2k_2)!(2k_3)!(2k_4)!}=-\frac{1}{120}k(k-1)(2k-3)(2k-1)(2k+1)\frac{\beta_{2k}}{(2k)!}\\
&\qquad\qquad -\frac{1}{12}(k-1)(2k-3)\frac{\beta_{2k-2}}{(2k-2)!}+\frac{1}{480}(2k-5)\frac{\beta_{2k-4}}{(2k-4)!},\\
&\sum\nolimits^{(4)} k_1^3\frac{\beta_{2k_1}\beta_{2k_2}\beta_{2k_3}\beta_{2k_4}}{(2k_1)!(2k_2)!(2k_3)!(2k_4)!}=-\frac{1}{240}k(k-1)(2k-3)(2k-1)(4k^2-6k+1)\frac{\beta_{2k}}{(2k)!}\\
&\qquad\qquad\qquad -\frac{1}{96}(2k-3)(6k^3-21k^2+17k+6)\frac{\beta_{2k-2}}{(2k-2)!}\\
&\qquad\qquad\qquad -\frac{1}{960}(2k-5)(13k-21)\frac{\beta_{2k-4}}{(2k-4)!},\\
&\sum\nolimits^{(4)} k_1^2k_2\frac{\beta_{2k_1}\beta_{2k_2}\beta_{2k_3}\beta_{2k_4}}{(2k_1)!(2k_2)!(2k_3)!(2k_4)!}=-\frac{1}{720}k(k-1)(2k-3)(2k-1)^2(2k+1)\frac{\beta_{2k}}{(2k)!}\\
&\qquad -\frac{1}{288}(k-1)(2k-3)(2k^2-k+6)\frac{\beta_{2k-2}}{(2k-2)!}+\frac{7}{2880}(k-3)(2k-5)\frac{\beta_{2k-4}}{(2k-4)!},\\
&\sum\nolimits^{(4)} k_1k_2k_3\frac{\beta_{2k_1}\beta_{2k_2}\beta_{2k_3}\beta_{2k_4}}{(2k_1)!(2k_2)!(2k_3)!(2k_4)!}=-\frac{1}{720}k(k-1)(k+1)(2k-3)(2k-1)(2k+1)\frac{\beta_{2k}}{(2k)!}\\
&\qquad +\frac{1}{288}(k-6)(k-1)(2k-3)(2k-1)\frac{\beta_{2k-2}}{(2k-2)!}-\frac{1}{2880}(2k-5)(4k-21)\frac{\beta_{2k-4}}{(2k-4)!}.
\end{align*}

%%%%%%%%%%%%%%%%%%%%%%%%%%%%%%%%%%%%%%%%%%%%%%%%%%%%%%%%%%%%%%%%%%%%%%%%%%%%%%%%%%%%%%%%%%%%%%%%%%%%%%%%%%%%%%%%%%%%

\subsection{Weighted sum formulas of $t$-values}

If $n=2$, we have
\begin{align*}
&\sum\nolimits^{(2)} t(2k_1)t(2k_2)=\frac{1}{2}(2k-1)t(2k),\\
&\sum\nolimits^{(2)} k_1t(2k_1)t(2k_2)=\frac{1}{4}k(2k-1)t(2k),\\
&\sum\nolimits^{(2)} k_1^2t(2k_1)t(2k_2)=\frac{1}{24}k(2k-1)(4k-1)t(2k)-\frac{1}{8}(2k-3)\zeta(2)t(2k-2),\\
&\sum\nolimits^{(2)} k_1k_2t(2k_1)t(2k_2)=\frac{1}{24}k(2k-1)(2k+1)t(2k)+\frac{1}{8}(2k-3)\zeta(2)t(2k-2),\\
&\sum\nolimits^{(2)} k_1^3t(2k_1)t(2k_2)=\frac{1}{16}k^2(2k-1)^2t(2k)-\frac{3}{16}k(2k-3)\zeta(2)t(2k-2),\\
&\sum\nolimits^{(2)} k_1^2k_2t(2k_1)t(2k_2)=\frac{1}{48}k^2(2k-1)(2k+1)t(2k)+\frac{1}{16}k(2k-3)\zeta(2)t(2k-2),\\
&\sum\nolimits^{(2)} k_1^4t(2k_1)t(2k_2)=\frac{1}{480}k(2k-1)(4k-1)(12k^2-6k-1)t(2k)\\
&\qquad\qquad -\frac{1}{32}(2k-3)(8k^2-6k+5)\zeta(2)t(2k-2)-\frac{3}{32}(2k-5)\zeta(4)t(2k-4),\\
&\sum\nolimits^{(2)} k_1^3k_2t(2k_1)t(2k_2)=\frac{1}{480}k(2k-1)(2k+1)(6k^2-1)t(2k)\\
&\qquad\qquad +\frac{1}{32}(2k-3)(2k^2-6k+5)\zeta(2)t(2k-2)+\frac{3}{32}(2k-5)\zeta(4)t(2k-4),\\
&\sum\nolimits^{(2)} k_1^2k_2^2t(2k_1)t(2k_2)=\frac{1}{480}k(2k-1)(2k+1)(4k^2+1)t(2k)\\
&\qquad\qquad +\frac{1}{32}(2k-3)(6k-5)\zeta(2)t(2k-2)-\frac{3}{32}(2k-5)\zeta(4)t(2k-4),\\
&\sum\nolimits^{(2)} k_1^5t(2k_1)t(2k_2)=\frac{1}{192}k^2(2k-1)^2(8k^2-4k-1)t(2k)\\
&\qquad\qquad -\frac{5}{64}k(2k-3)(4k^2-6k+5)\zeta(2)t(2k-2)-\frac{15}{64}k(2k-5)\zeta(4)t(2k-4),\\
&\sum\nolimits^{(2)} k_1^4k_2t(2k_1)t(2k_2)=\frac{1}{960}k^2(2k-1)(2k+1)(8k^2-3)t(2k)\\
&\qquad\qquad +\frac{1}{64}k(2k-3)(4k^2-18k+15)\zeta(2)t(2k-2)+\frac{9}{64}k(2k-5)\zeta(4)t(2k-4),\\
&\sum\nolimits^{(2)} k_1^3k_2^2t(2k_1)t(2k_2)=\frac{1}{960}k^2(2k-1)(2k+1)(4k^2+1)t(2k)\\
&\qquad\qquad +\frac{1}{64}k(2k-3)(6k-5)\zeta(2)t(2k-2)-\frac{3}{64}k(2k-5)\zeta(4)t(2k-4).
\end{align*}

If $n=3$, we have
\begin{align*}
&\sum\nolimits^{(3)} t(2k_1)t(2k_2)t(2k_3)=\frac{1}{4}(k-1)(2k-1)t(2k)-\frac{3}{8}\zeta(2)t(2k-2),\\
&\sum\nolimits^{(3)} k_1t(2k_1)t(2k_2)t(2k_3)=\frac{1}{12}k(k-1)(2k-1)t(2k)-\frac{1}{8}k\zeta(2)t(2k-2),\\
&\sum\nolimits^{(3)} k_1^2t(2k_1)t(2k_2)t(2k_3)=\frac{1}{48}k(k-1)(2k-1)^2t(2k)\\
&\qquad\qquad\qquad\qquad\qquad\qquad -\frac{1}{16}(4k^2-11k+9)\zeta(2)t(2k-2),\\
&\sum\nolimits^{(3)} k_1k_2t(2k_1)t(2k_2)t(2k_3)=\frac{1}{96}k(k-1)(2k-1)(2k+1)t(2k)\\
&\qquad\qquad\qquad\qquad\qquad\qquad +\frac{1}{32}(k-1)(2k-9)\zeta(2)t(2k-2),\\
&\sum\nolimits^{(3)} k_1^3t(2k_1)t(2k_2)t(2k_3)=\frac{1}{480}k(k-1)(2k-1)(12k^2-12k+1)t(2k)\\
&\qquad\qquad -\frac{1}{32}(8k^3-24k^2+17k+3)\zeta(2)t(2k-2)+\frac{3}{32}(2k-5)\zeta(4)t(2k-4),\\
&\sum\nolimits^{(3)} k_1^2k_2t(2k_1)t(2k_2)t(2k_3)=\frac{1}{960}k(k-1)(2k-1)(2k+1)(4k-1)t(2k)\\
&\qquad\qquad -\frac{1}{64}(k-1)(2k+3)\zeta(2)t(2k-2)-\frac{3}{64}(2k-5)\zeta(4)t(2k-4),\\
&\sum\nolimits^{(3)} k_1k_2k_3t(2k_1)t(2k_2)t(2k_3)=\frac{1}{480}k(k-1)(k+1)(2k-1)(2k+1)t(2k)\\
&\qquad\qquad +\frac{1}{32}(k-3)(k-1)(2k-1)\zeta(2)t(2k-2)+\frac{3}{32}(2k-5)\zeta(4)t(2k-4),\\
&\sum\nolimits^{(3)} k_1^4t(2k_1)t(2k_2)t(2k_3)=\frac{1}{960}k(k-1)(2k-1)^2(8k^2-8k-1)t(2k)\\
&\quad -\frac{1}{64}(16k^4-64k^3+96k^2-79k+39)\zeta(2)t(2k-2)+\frac{3}{64}(2k-5)(3k+1)\zeta(4)t(2k-4),\\
&\sum\nolimits^{(3)} k_1^3k_2t(2k_1)t(2k_2)t(2k_3)=\frac{1}{1920}k(k-1)(2k-1)(2k+1)(4k^2-2k-1)t(2k)\\
&\quad -\frac{1}{128}(k-1)(16k^2-46k+39)\zeta(2)t(2k-2)-\frac{3}{128}(k+1)(2k-5)\zeta(4)t(2k-4),\\
&\sum\nolimits^{(3)} k_1^2k_2^2t(2k_1)t(2k_2)t(2k_3)=\frac{1}{5760}k(k-1)(2k-1)(2k+1)(8k^2-4k+3)t(2k)\\
&\qquad\qquad\qquad -\frac{1}{384}(k-1)(8k^3-64k^2+168k-117)\zeta(2)t(2k-2)\\
&\qquad\qquad\qquad-\frac{1}{128}(2k-5)(7k-3)\zeta(4)t(2k-4),\\
&\sum\nolimits^{(3)} k_1^2k_2k_3t(2k_1)t(2k_2)t(2k_3)=\frac{1}{1440}k^2(k-1)(k+1)(2k-1)(2k+1)t(2k)\\
&\qquad +\frac{1}{96}k(k-3)(k-1)(2k-1)\zeta(2)t(2k-2)+\frac{1}{32}k(2k-5)\zeta(4)t(2k-4).
\end{align*}

If $n=4$, we have
\begin{align*}
&\sum\nolimits^{(4)} t(2k_1)t(2k_2)t(2k_3)t(2k_4)=\frac{1}{24}(k-1)(2k-3)(2k-1)t(2k)\\
&\qquad\qquad\qquad\qquad\qquad-\frac{1}{4}(2k-3)\zeta(2)t(2k-2),\\
&\sum\nolimits^{(4)} k_1t(2k_1)t(2k_2)t(2k_3)t(2k_4)=\frac{1}{96}k(k-1)(2k-3)(2k-1)t(2k)\\
&\qquad\qquad\qquad\qquad\qquad-\frac{1}{16}k(2k-3)\zeta(2)t(2k-2),\\
&\sum\nolimits^{(4)} k_1^2t(2k_1)t(2k_2)t(2k_3)t(2k_4)=\frac{1}{960}k(k-1)(2k-3)(2k-1)(4k-3)t(2k)\\
&\qquad\qquad-\frac{1}{16}(2k-3)(k^2-3k+3)\zeta(2)t(2k-2)+\frac{9}{128}(2k-5)\zeta(4)t(2k-4),\\
&\sum\nolimits^{(4)} k_1k_2t(2k_1)t(2k_2)t(2k_3)t(2k_4)=\frac{1}{960}k(k-1)(2k-3)(2k-1)(2k+1)t(2k)\\
&\qquad\qquad-\frac{1}{16}(k-1)(2k-3)\zeta(2)t(2k-2)-\frac{3}{128}(2k-5)\zeta(4)t(2k-4),\\
&\sum\nolimits^{(4)} k_1^3t(2k_1)t(2k_2)t(2k_3)t(2k_4)=\frac{1}{1920}k(k-1)(2k-3)(2k-1)(4k^2-6k+1)t(2k)\\
&\qquad\qquad\qquad-\frac{1}{128}(2k-3)(6k^3-21k^2+17k+6)\zeta(2)t(2k-2)\\
&\qquad\qquad\qquad+\frac{3}{256}(2k-5)(13k-21)\zeta(4)t(2k-4),\\
&\sum\nolimits^{(4)} k_1^2k_2t(2k_1)t(2k_2)t(2k_3)t(2k_4)=\frac{1}{5760}k(k-1)(2k-3)(2k-1)^2(2k+1)t(2k)\\
&\quad-\frac{1}{384}(k-1)(2k-3)(2k^2-k+6)\zeta(2)t(2k-2)-\frac{7}{256}(k-3)(2k-5)\zeta(4)t(2k-4),\\
&\sum\nolimits^{(4)} k_1k_2k_3t(2k_1)t(2k_2)t(2k_3)t(2k_4)=\frac{1}{5760}k(k-1)(k+1)(2k-3)(2k-1)(2k+1)t(2k)\\
&\qquad\qquad\qquad+\frac{1}{384}(k-6)(k-1)(2k-3)(2k-1)\zeta(2)t(2k-2)\\
&\qquad\qquad\qquad+\frac{1}{256}(2k-5)(4k-21)\zeta(4)t(2k-4).
\end{align*}

%%--------------------------------------------------------------------------------------------------------------
%%%%%%%%%%%%%%%%%%%%%%%%%%%%%%%%%%%%%%%%%%%%%%%%%%%%%%%%%%%%%%%%%%%%%%%%%%%%%%%%%%%%%%%%%%%%%%%%%%%%%%%%%%%%%%%%%%%%

\subsection{Weighted sum formulas of multiple $t$-values}

If $n=2$, we have
\begin{align*}
&\sum\nolimits^{(2)} t(2k_1,2k_2)=\frac{1}{4}t(2k),\\
&\sum\nolimits^{(2)}(k_1^2+k_2^2) t(2k_1,2k_2)=\frac{1}{8}k(2k-1)t(2k)-\frac{1}{8}(2k-3)\zeta(2)t(2k-2),\\
&\sum\nolimits^{(2)} k_1k_2 t(2k_1,2k_2)=\frac{1}{16}kt(2k)+\frac{1}{16}(2k-3)\zeta(2)t(2k-2),\\
&\sum\nolimits^{(2)}(k_1^3+k_2^3) t(2k_1,2k_2)=\frac{1}{16}k^2(4k-3)t(2k)-\frac{3}{16}k(2k-3)\zeta(2)t(2k-2),\\
&\sum\nolimits^{(2)}(k_1^4+k_2^4) t(2k_1,2k_2)=\frac{1}{32}k(2k-1)(4k^2-2k-1)t(2k)\\
&\qquad\qquad -\frac{1}{32}(2k-3)(8k^2-6k+5)\zeta(2)t(2k-2)-\frac{3}{32}(2k-5)\zeta(4)t(2k-4),\\
&\sum\nolimits^{(2)}(k_1^3k_2+k_1k_2^3) t(2k_1,2k_2)=\frac{1}{32}k(2k^2-1)t(2k)\\
&\qquad\qquad +\frac{1}{32}(2k-3)(2k^2-6k+5)\zeta(2)t(2k-2)+\frac{3}{32}(2k-5)\zeta(4)t(2k-4),\\
&\sum\nolimits^{(2)} k_1^2k_2^2 t(2k_1,2k_2)=\frac{1}{64}kt(2k)+\frac{1}{64}(2k-3)(6k-5)\zeta(2)t(2k-2)\\
&\qquad\qquad -\frac{3}{64}(2k-5)\zeta(4)t(2k-4),\\
&\sum\nolimits^{(2)}(k_1^5+k_2^5) t(2k_1,2k_2)=\frac{1}{64}k^2(16k^3-20k^2+5)t(2k)\\
&\qquad\qquad -\frac{5}{64}k(2k-3)(4k^2-6k+5)\zeta(2)t(2k-2)-\frac{15}{64}k(2k-5)\zeta(4)t(2k-4),\\
&\sum\nolimits^{(2)}(k_1^4k_2+k_1k_2^4) t(2k_1,2k_2)=\frac{1}{64}k^2(4k^2-3)t(2k)\\
&\qquad\qquad +\frac{1}{64}k(2k-3)(4k^2-18k+15)\zeta(2)t(2k-2)+\frac{9}{64}k(2k-5)\zeta(4)t(2k-4)
\end{align*}
and
\begin{align*}
&\sum\nolimits^{(2)} t^{\star}(2k_1,2k_2)=\frac{1}{4}(4k-3)t(2k),\\
&\sum\nolimits^{(2)}(k_1^2+k_2^2) t^{\star}(2k_1,2k_2)=\frac{1}{24}k(2k-1)(8k-5)t(2k)-\frac{1}{8}(2k-3)\zeta(2)t(2k-2),\\
&\sum\nolimits^{(2)} k_1k_2 t^{\star}(2k_1,2k_2)=\frac{1}{48}k(8k^2-5)t(2k)+\frac{1}{16}(2k-3)\zeta(2)t(2k-2),\\
&\sum\nolimits^{(2)}(k_1^3+k_2^3) t^{\star}(2k_1,2k_2)=\frac{1}{16}k^2(8k^2-12k+5)t(2k)-\frac{3}{16}k(2k-3)\zeta(2)t(2k-2),\\
&\sum\nolimits^{(2)}(k_1^4+k_2^4) t^{\star}(2k_1,2k_2)=\frac{1}{480}k(2k-1)(96k^3-132k^2+34k+17)t(2k)\\
&\qquad\qquad -\frac{1}{32}(2k-3)(8k^2-6k+5)\zeta(2)t(2k-2)-\frac{3}{32}(2k-5)\zeta(4)t(2k-4),\\
&\sum\nolimits^{(2)}(k_1^3k_2+k_1k_2^3) t^{\star}(2k_1,2k_2)=\frac{1}{480}k(48k^4-50k^2+17)t(2k)\\
&\qquad\qquad +\frac{1}{32}(2k-3)(2k^2-6k+5)\zeta(2)t(2k-2)+\frac{3}{32}(2k-5)\zeta(4)t(2k-4),\\
&\sum\nolimits^{(2)} k_1^2k_2^2 t^{\star}(2k_1,2k_2)=\frac{1}{960}k(32k^4-17)t(2k)+\frac{1}{64}(2k-3)(6k-5)\zeta(2)t(2k-2)\\
&\qquad\qquad\qquad\qquad -\frac{3}{64}(2k-5)\zeta(4)t(2k-4),\\
&\sum\nolimits^{(2)}(k_1^5+k_2^5) t^{\star}(2k_1,2k_2)=\frac{1}{192}k^2(64k^4-144k^3+100k^2-17)t(2k)\\
&\qquad\qquad -\frac{5}{64}k(2k-3)(4k^2-6k+5)\zeta(2)t(2k-2)-\frac{15}{64}k(2k-5)\zeta(4)t(2k-4),\\
&\sum\nolimits^{(2)}(k_1^4k_2+k_1k_2^4) t^{\star}(2k_1,2k_2)=\frac{1}{960}k^2(64k^4-100k^2+51)t(2k)\\
&\qquad\qquad +\frac{1}{64}k(2k-3)(4k^2-18k+15)\zeta(2)t(2k-2)+\frac{9}{64}k(2k-5)\zeta(4)t(2k-4).
\end{align*}

If $n=3$, we have
\begin{align*}
&\sum\nolimits^{(3)} t(2k_1,2k_2,2k_3)=\frac{1}{8}t(2k)-\frac{1}{16}\zeta(2)t(2k-2),\\
&\sum\nolimits^{(3)} (k_1^2+k_2^2+k_3^2)t(2k_1,2k_2,2k_3)=\frac{1}{32}k(4k-3)t(2k)-\frac{1}{32}(2k^2-3)\zeta(2)t(2k-2),\\
&\sum\nolimits^{(3)} (k_1k_2+k_1k_3+k_2k_3)t(2k_1,2k_2,2k_3)=\frac{3}{64}kt(2k)-\frac{3}{64}\zeta(2)t(2k-2),\\
&\sum\nolimits^{(3)} (k_1^3+k_2^3+k_3^3)t(2k_1,2k_2,2k_3)=\frac{1}{128}k(16k^2-18k+3)t(2k)\\
&\qquad\qquad-\frac{1}{128}(8k^3-18k+3)\zeta(2)t(2k-2)+\frac{3}{128}(2k-5)\zeta(4)t(2k-4),\\
&\sum\nolimits^{(3)} \sum\limits_{1\leqslant i<j\leqslant 3}(k_i^2k_j+k_ik_j^2)t(2k_1,2k_2,2k_3)=\frac{3}{128}k(2k-1)t(2k)\\
&\qquad\qquad-\frac{3}{128}(2k-1)\zeta(2)t(2k-2)-\frac{3}{128}(2k-5)\zeta(4)t(2k-4),\\
&\sum\nolimits^{(3)} k_1k_2k_3t(2k_1,2k_2,2k_3)=\frac{1}{128}kt(2k)-\frac{1}{128}\zeta(2)t(2k-2)\\
&\qquad\qquad\qquad\qquad\qquad +\frac{1}{128}(2k-5)\zeta(4)t(2k-4),\\
&\sum\nolimits^{(3)} (k_1^4+k_2^4+k_3^4)t(2k_1,2k_2,2k_3)=\frac{1}{128}k(16k^3-24k^2+6k+3)t(2k)\\
&\qquad-\frac{1}{128}(8k^4-36k^2+42k-21)\zeta(2)t(2k-2)+\frac{3}{128}(2k-5)(2k-3)\zeta(4)t(2k-4),\\
&\sum\nolimits^{(3)} \sum\limits_{1\leqslant i<j\leqslant 3}(k_i^3k_j+k_ik_j^3)t(2k_1,2k_2,2k_3)=\frac{3}{128}k(k-1)(2k+1)t(2k)\\
&\qquad\qquad-\frac{3}{128}(k-1)(6k-7)\zeta(2)t(2k-2)-\frac{3}{128}(k-3)(2k-5)\zeta(4)t(2k-4),\\
&\sum\nolimits^{(3)} (k_1^2k_2^2+k_1^2k_3^2+k_2^2k_3^2)t(2k_1,2k_2,2k_3)=-\frac{1}{256}k(2k-3)t(2k)\\
&\qquad+\frac{1}{256}(12k^2-34k+21)\zeta(2)t(2k-2)-\frac{1}{256}(2k-5)(2k+9)\zeta(4)t(2k-4)
\end{align*}
and
\begin{align*}
&\sum\nolimits^{(3)} t^{\star}(2k_1,2k_2,2k_3)=\frac{1}{8}(4k^2-10k+5)t(2k)-\frac{1}{16}\zeta(2)t(2k-2),\\
&\sum\nolimits^{(3)} (k_1^2+k_2^2+k_3^2)t^{\star}(2k_1,2k_2,2k_3)=\frac{1}{96}k(24k^3-80k^2+78k-25)t(2k)\\
&\qquad\qquad\qquad\qquad -\frac{1}{32}(6k^2-22k+21)\zeta(2)t(2k-2),\\
&\sum\nolimits^{(3)} (k_1k_2+k_1k_3+k_2k_3)t^{\star}(2k_1,2k_2,2k_3)=\frac{1}{192}k(24k^3-40k^2-18k+25)t(2k)\\
&\qquad\qquad\qquad\qquad+\frac{1}{64}(4k^2-22k+21)\zeta(2)t(2k-2),\\
&\sum\nolimits^{(3)} (k_1^3+k_2^3+k_3^3)t^{\star}(2k_1,2k_2,2k_3)=\frac{1}{640}k(6k-1)(16k^3-64k^2+76k-29)t(2k)\\
&\qquad\qquad-\frac{1}{128}(24k^3-96k^2+86k+9)\zeta(2)t(2k-2)+\frac{9}{128}(2k-5)\zeta(4)t(2k-4),\\
&\sum\nolimits^{(3)} \sum\limits_{1\leqslant i<j\leqslant 3}(k_i^2k_j+k_ik_j^2)t^{\star}(2k_1,2k_2,2k_3)=\frac{1}{1920}k(2k-1)(96k^3-152k^2-76k+87)t(2k)\\
&\qquad\qquad-\frac{1}{128}(8k^2-2k-9)\zeta(2)t(2k-2)-\frac{9}{128}(2k-5)\zeta(4)t(2k-4),\\
&\sum\nolimits^{(3)} k_1k_2k_3t^{\star}(2k_1,2k_2,2k_3)=\frac{1}{1920}k(16k^4-60k^2+29)t(2k)\\
&\qquad\qquad +\frac{1}{384}(8k^3-36k^2+40k-9)\zeta(2)t(2k-2)+\frac{3}{128}(2k-5)\zeta(4)t(2k-4),\\
&\sum\nolimits^{(3)} (k_1^4+k_2^4+k_3^4)t^{\star}(2k_1,2k_2,2k_3)=\frac{1}{1920}k(192k^5-960k^4+1560k^3-1000k^2+108k\\
&\qquad\qquad +85)t(2k)-\frac{1}{128}(24k^4-128k^3+228k^2-200k+99)\zeta(2)t(2k-2)\\
&\qquad\qquad+\frac{3}{128}(2k-5)(4k+5)\zeta(4)t(2k-4),\\
&\sum\nolimits^{(3)} \sum\limits_{1\leqslant i<j\leqslant 3}(k_i^3k_j+k_ik_j^3)t^{\star}(2k_1,2k_2,2k_3)=\frac{1}{1920}k(k-1)(96k^4-144k^3-144k^2\\
&\qquad\qquad +106k+85)t(2k)-\frac{1}{128}(k-1)(32k^2-110k+99)\zeta(2)t(2k-2)\\
&\qquad\qquad-\frac{3}{128}(k+5)(2k-5)\zeta(4)t(2k-4),\\
&\sum\nolimits^{(3)} (k_1^2k_2^2+k_1^2k_3^2+k_2^2k_3^2)t(2k_1,2k_2,2k_3)=\frac{1}{3840}k(64k^5-160k^4+120k^3-124k\\
&\qquad\qquad +85)t(2k)-\frac{1}{768}(16k^4-144k^3+500k^2-672k+297)\zeta(2)t(2k-2)\\
&\qquad\qquad -\frac{3}{256}(2k-5)(4k-5)\zeta(4)t(2k-4).
\end{align*}

If $n=4$, we have
\begin{align*}
&\sum\nolimits^{(4)} t(2k_1,2k_2,2k_3,2k_4)=\frac{5}{64}t(2k)-\frac{3}{64}\zeta(2)t(2k-2),\\
&\sum\nolimits^{(4)} (k_1^2+k_2^2+k_3^2+k_4^2)t(2k_1,2k_2,2k_3,2k_4)=\frac{1}{128}k(10k-9)t(2k)\\
&\qquad\qquad-\frac{3}{128}(2k^2-k-2)\zeta(2)t(2k-2)+\frac{3}{256}(2k-5)\zeta(4)t(2k-4),\\
&\sum\nolimits^{(4)} \sum\limits_{1\leqslant i<j\leqslant 4}k_ik_j t(2k_1,2k_2,2k_3,2k_4)=\frac{9}{256}kt(2k)\\
&\qquad\qquad-\frac{3}{256}(k+2)\zeta(2)t(2k-2)-\frac{3}{512}(2k-5)\zeta(4)t(2k-4),\\
&\sum\nolimits^{(4)} (k_1^3+k_2^3+k_3^3+k_4^3)t(2k_1,2k_2,2k_3,2k_4)=\frac{1}{512}k(40k^2-54k+15)t(2k)\\
&\qquad\qquad-\frac{3}{512}(8k^3-6k^2-10k+3)\zeta(2)t(2k-2)+\frac{9}{512}k(2k-5)\zeta(4)t(2k-4),\\
&\sum\nolimits^{(4)} \sum\limits_{1\leqslant i<j\leqslant 4}(k_i^2k_j+k_ik_j^2) t(2k_1,2k_2,2k_3,2k_4)=\frac{3}{512}k(6k-5)t(2k)\\
&\qquad\qquad-\frac{3}{512}(2k^2+2k-3)\zeta(2)t(2k-2)-\frac{3}{512}k(2k-5)\zeta(4)t(2k-4),\\
&\sum\nolimits^{(4)} (k_1k_2k_3+k_1k_2k_4+k_1k_3k_4+k_2k_3k_4) t(2k_1,2k_2,2k_3,2k_4)=\frac{5}{512}kt(2k)\\
&\qquad\qquad\qquad\qquad\qquad -\frac{1}{512}(2k+3)\zeta(2)t(2k-2)
\end{align*}
and
\begin{align*}
&\sum\nolimits^{(4)} t^{\star}(2k_1,2k_2,2k_3,2k_4)=\frac{1}{192}(4k-7)(8k^2-28k+15)t(2k)\\
&\qquad\qquad\qquad\qquad\qquad -\frac{1}{64}(4k-9)\zeta(2)t(2k-2),\\
&\sum\nolimits^{(4)} (k_1^2+k_2^2+k_3^2+k_4^2)t^{\star}(2k_1,2k_2,2k_3,2k_4)=\frac{1}{1920}k(128k^4-840k^3+1880k^2\\
&\qquad\qquad -1680k+527)t(2k)-\frac{1}{384}(32k^3-210k^2+469k-348)\zeta(2)t(2k-2)\\
&\qquad\qquad-\frac{1}{256}(2k-5)\zeta(4)t(2k-4),\\
&\sum\nolimits^{(4)} \sum\limits_{1\leqslant i<j\leqslant 4}k_ik_j t^{\star}(2k_1,2k_2,2k_3,2k_4)=\frac{1}{3840}k(192k^4-840k^3+680k^2+630k\\
&\qquad\qquad -527)t(2k)+\frac{1}{768}(8k^3-156k^2+469k-348)\zeta(2)t(2k-2)\\
&\qquad\qquad +\frac{1}{512}(2k-5)\zeta(4)t(2k-4),\\
&\sum\nolimits^{(4)} (k_1^3+k_2^3+k_3^3+k_4^3)t^{\star}(2k_1,2k_2,2k_3,2k_4)=\frac{1}{7680}k(256k^5-2016k^4+5640k^3\\
&\qquad\qquad -6720k^2+3464k-609)t(2k)-\frac{1}{512}(32k^4-232k^3+534k^2-344k\\
&\qquad\qquad -69)\zeta(2)t(2k-2)+\frac{3}{512}(2k-5)(5k-26)\zeta(4)t(2k-4),\\
&\sum\nolimits^{(4)} \sum\limits_{1\leqslant i<j\leqslant 4}(k_i^2k_j+k_ik_j^2) t^{\star}(2k_1,2k_2,2k_3,2k_4)=\frac{1}{7680}k(256k^5-1344k^4+1880k^3\\
&\qquad\qquad -1356k+609)t(2k)-\frac{1}{1536}(32k^4-144k^3+274k^2-360k\\
&\qquad\qquad +207)\zeta(2)t(2k-2)-\frac{1}{512}(2k-5)(17k-78)\zeta(4)t(2k-4),\\
&\sum\nolimits^{(4)} (k_1k_2k_3+k_1k_2k_4+k_1k_3k_4+k_2k_3k_4) t^{\star}(2k_1,2k_2,2k_3,2k_4)=\frac{1}{23040}k(128k^5-336k^4\\
&\qquad\qquad -520k^3+1260k^2+302k-609)t(2k)+\frac{1}{1536}(16k^4-152k^3+404k^2-352k\\
&\qquad\qquad +69)\zeta(2)t(2k-2)+\frac{1}{256}(2k-5)(3k-13)\zeta(4)t(2k-4).
\end{align*}

\noindent {\bf Acknowledgments.}  The authors express their deep gratitude to Professor Masanobu Kaneko for valuable discussions and comments. The second author is supported by the China Scholarship Council (No. 201806310063).


{\small
\begin{thebibliography}{99}

\bibitem{ELO2017}
M. Eie, W-C. Liaw and Y.L. Ong, Several weighted sum formula of multiple zeta values, Int. J. Number Theory, 2017, {\bf 13}(9), 2253-2264.

\bibitem{ELO2018}
M. Eie and Y.L. Ong, Sum formulas of multiple zeta values with even arguments and polynomial weights, J. Number Theory, 2018, {\bf 188}, 247-262.

\bibitem{GKZ2006}
H. Gangl, M. Kaneko and D. Zagier, Double zeta values and modular forms, in: Proceedings of the Conference in Memory of Tsuneo Arakawa, 2006, 71-106.

\bibitem{G2016}
M. Gen${\check{\rm c}}$ev, On restricted sum formulas for multiple zeta values with even arguments, Arch. Math., 2016, {\bf 107}, 9-22.

\bibitem{GL2015}
L. Guo, P. Lei and J. Zhao, Families of weighted sum formulas for multiple zeta values, Int. J. Number Theory, 2015, {\bf 11}(3), 997-1025.

\bibitem{H1992}
M.E. Hoffman, Multiple harmonic series, Pacific J. Math., 1992, {\bf 152}(2), 275-290.

\bibitem{H2016}
M.E. Hoffman, An odd variant of multiple zeta values, Comm. Number Theory Phys., to appear, arXiv:1612.05232v4 [math.NT].

\bibitem{H2017}
M.E. Hoffman, On multiple zeta values of even arguments, Int. J. Number Theory, 2017, {\bf 13}(3), 705-716.

\bibitem{KMT2014}
Y. Komori, K. Matsumoto and H. Tsumura, A study on multiple zeta values from viewpoint of zeta functions of root systems, Funct. Approx. Comment. Math., 2015, {\bf 51}, 47-76.

\bibitem{LQ2016} Z. Li and C. Qin, Some relations deduced from regularized double shuffle relations of multiple zeta values, preprint, arXiv: 1610.05480 [math.NT].

\bibitem{LQ2019}
Z. Li and C. Qin, Weighted sum formulas of multiple zeta values with even arguments, Math. Z., 2019, {\bf 291}, 1337-1356.

\bibitem{N2009}
T. Nakamura, Restricted and weighted sum formulas for double zeta values of even weight, $\check{\rm S}$iauliai Math. Semin., 2009, {\bf 12}, 151-155.

\bibitem{SC2011}
Z. Shen and T. Cai, Some identities for multiple Hurwitz zeta values (in Chinese), Sci. Sinica Math., 2011, {\bf 41}, 955-970.

\bibitem{SC2012}
Z. Shen and T. Cai, Some identities for multiple zeta values, J. Number Theory, 2012, {\bf 132}, 314-323.

\bibitem{SJ2017}
Z. Shen and L. Jia, Some identities for multiple Hurwitz zeta values, J. Number Theory, 2017, {\bf 179}, 256-267.

\bibitem{DZ1994}
D. Zagier, Values of zeta functions and their applications, First European Congress
of Mathematics, Volume II, Birkhauser, Boston, 1994, {\bf 120}, 497-512.

\bibitem{Z2015}
J. Zhao, Sum formula of multiple Hurwitz-zeta values, Forum. Math., 2015, {\bf 27}, 929-936.

\bibitem{Z2016}
J. Zhao, Multiple zeta functions, multiple polylogarithms and their special values, Series on Number
Theory and its Applications, 12, World Scientific Publishing Co. Pte. Ltd., Hackensack, NJ, 2016.


\end{thebibliography}}
 \end{document}